\documentclass{article}
\usepackage{amssymb}
\usepackage{amsmath}

\title{Semi-stable extensions on arithmetic surfaces}
\author{C. Soul\'e }
\date{ \ }

\begin{document}

\maketitle

Let $S$ be a smooth  projective curve over the complex numbers and $X \to S$ a semi-stable projective family of curves. Assume that both $S$ and the generic fiber of $X$ over $S$ have genus at least two. Then the sheaf of absolute differentials $\Omega_X^1$ defines a vector bundle on $X$ which is semi-stable in the sense of Mumford-Nakano with respect to the canonical line bundle on $X$. The Bogomolov inequality
$$
c_1^2 (\Omega_X^1) \leq 4 \, c_2 (\Omega_X^1)
$$
leads to an upper bound for the self-intersection $c_1 (\omega_{X/S})^2$ of the relative dualizing sheaf $\omega_{X/S}$.

\medskip

Assume now that $S$ is the spectrum ${\rm Spec} \, ({\mathcal O}_F)$ of the ring of integers in a number field $F$ and that $X \to S$ is a semi-stable (regular) curve over $S$, with generic genus at least two. In \cite{P}, Parshin asked for a similar upper bound for the arithmetic self-intersection $\hat c_1 (\bar\omega_{X/S})^2$ of the relative dualizing sheaf of $X$ over $S$, equipped with its Arakelov metric. He and Moret-Bailly \cite{MB} proved that a good upper bound for this real number $\hat c_1 (\bar\omega_{X/S})^2$ would have beautiful arithmetic consequences (including the $abc$ conjecture).

\medskip

If one tries to mimick in the arithmetic case the proof that we have just checked in the geometric case, one soon faces the difficulty that we do not know any arithmetic analog for the sheaf of absolute differentials $\Omega_X^1$. In \cite{M1}, Miyaoka proposed to turn this difficulty as follows. He noticed that, in the geometric case, any general enough rank two extension $E$ of $\omega_{X/S}$ by the pull-back to $X$ of $\Omega_S^1$ is semi-stable and that it can be used instead of $\Omega_X^1$ in the argument. When $S = {\rm Spec} \, ({\mathcal O}_F)$ it is then natural to apply an arithmetic analog of Bogomolov inequality to a rank two extension $\bar E$ of $\bar\omega_{X/S}$ by some hermitian line bundle pulled back from $S$. 

\medskip

But then, a new difficulty arises. Namely, the second Chern number $\hat c_2 (\bar E)$ of $\bar E$ is more involved in the arithmetic case than in the geometric one, as it contains an archimedean summand -- an integral over the set of complex points of $X$ -- which is not easy to bound from above.

\medskip

In this paper, which is a sequel to \cite{S1} and \cite{S2},  although we are unable to prove Parshin's conjecture by Miyaoka's argument,
we show that his method still provides interesting lower bounds  form
 some successive minima of the euclidean lattice of sections of hermitian line bundles on the arithmetic surface $X$.

\medskip

More precisely, we consider an hermitian line bundle $\bar N$ on $X$, with positive degree on the generic fiber and such that $N \otimes \omega_{X/S}$ is the square of a line bundle on $X$. We prove that, when $k$ is big enough, the logarithm $\mu_k$ of the $k$-th successive minimum of $H^1 (X , N^{-1})$, endowed with its $L^2$-metric, is bounded below:
$$
\mu_k \geq \frac{\hat c_1 (\bar N)^2}{2 \, n \, d} - A \, , \leqno (*)
$$
where $n$ is the degree of $N$, $d = [F: {\mathbb Q}]$ and $A$ is a simple constant (Theorem 2).

\medskip

This result is a complement to Theorem 4 in \cite{S2}, where smaller values of $k$ were studied. The proof is similar and consists mainly in making precise Miyaoka's assertion that a general extension $E$ of $N$ by the trivial line bundle is semi-stable on $X \otimes \bar F$. For that, again inspired by Miyaoka, we write $E$ as an extension
$$
0 \to L \to E \to M \to 0
$$
over $X \otimes \bar F$, with $L = \omega_{X/S} \otimes M^{-1}$ being the Serre dual of $M$, and we show (Proposition 1) that $E$ is semi-stable on $X \otimes \bar F$ as soon as the boundary map
$$
\partial : H^0 (X_{\bar F} , M) \to H^1 (X_{\bar F} , L)
$$
is an isomorphism. Next, we give an upper bound for the dimension of a vector space $V \subset {\rm Ext} \, (L,M)$ such that, for every extension class in $V$, the corresponding map $\partial$ is singular (Proposition 2). By a standard argument it follows that, if $k$ is big enough, there exists an extension $E$ of $N$ by ${\mathcal O}_X$ which is semi-stable over $X \otimes \bar F$ and such that the $L^2$-norm of its extension class is bounded above by $\exp \, (\mu_k + A)$. The proof of $(*)$ (see Theorem 2) then follows from a theorem ``\`a la Bogomolov'' for semi-stable hermitian vector bundles on arithmetic surfaces, which is due to Miyaoka \cite{M2}, \cite{S1} and Moriwaki \cite{Mo}.

\medskip

The geometric aspect of our argument can also be expressed in terms of the secant variety $\Sigma_d$ of a smooth projective curve $C$. In \cite{V}, Voisin gave an upper bound for the dimension of projective spaces contained in $\Sigma_d$, when $d$ is small enough with respect to the degree of $C$. In Theorem 1, we prove a similar result for a slightly bigger value of $d$.

\medskip

When doing this work, I got help from C.~Gasbarri, B.~Mazur, Y.~Miyaoka and especially C.~Voisin, who
found a gap in the proof of the key Proposition 2 and fixed it. I wish to
express to them my gratitude, as well as to the organizers of this conference.

\vglue 1cm

\noindent {\bf Notation.} Given two line bundles on a scheme $X$, we denote by $L^{-1}$ the dual of $L$ and by $LM$ the tensor product of $L$ with $M$.

\section{Semi-stable extensions on curves}

\subsection{ \ } 

Let $k = \bar k$ be an algebraically closed field, and $C$ a smooth connected projective curve of genus $g \geq 0$ over $k$. Let $L$ and $M$ be two line bundles on $X$ and
$$
0 \to L \to E \to M \to 0 \leqno (1)
$$
a rank two extension of $M$ by $L$. Consider the associated boundary map in cohomology
$$
\partial : H^0 (C,M) \to H^1 (C,L) \, .
$$

\medskip

\noindent {\bf Proposition 1.} Assume that $\deg (L) \leq \deg (M)$ and that
\begin{enumerate}
\item[a)] Either $\deg (L) + \deg (M) \geq 2g-2$ and $\partial$ is injective;
\item[b)] Or $\deg (L) + \deg (M) \leq 2g-2$ and $\partial$ is surjective.
\end{enumerate}
Then the vector bundle $E$ is semi-stable on $C$.

\subsection{Proof.} Let us prove a) by contradiction. Let $N \subset E$ be a line bundle on $C$ such that
$$
\deg (N) > \frac{\deg (E)}{2} = \frac{\deg (L) + \deg (M)}{2} \, .
$$
Then $\deg (N) > \deg (L)$, therefore $N \cap L = 0$, and the composite map
$$
N \to E \to M
$$
is injective. The extension
$$
0 \to L \to E' \to N \to 0
$$
induced by (1) and this map is split. Therefore the associated boundary map
$$
H^0 (C,N) \to H^1 (C,L)
$$
is zero, {\it i.e.} the restriction of $\partial$ to $H^0 (C,N) \subset H^0 (C,M)$ vanishes.

\medskip

On the other hand, since $\deg (L) + \deg (M) \geq 2g-2$, we have
$$
\deg (N) > g-1 \, ,
$$
hence, by Riemann-Roch, $H^0 (C,N) \ne 0$. This contradicts the assumption that $\partial$ is injective.

\medskip

To prove b) by contradiction we may consider a quotient $N$ of $E$ of degree less than $\deg (E) / 2$ and look at the extension
$$
0 \to N \to E' \to M \to 0
$$
induced by the composite map $L \to E' \to N$. Alternatively, one can deduce b) from a) by considering the Serre dual of $E$.

\subsection{Remark.} There are cases where $E$ is semi-stable when neither a) nor b) holds.

\subsubsection{ \ } For instance, when $C$ is an elliptic curve and $A \in C(k)$, let $L = {\mathcal O} (-A)$ and $M = {\mathcal O} (A)$. The group of extensions
$$
{\rm Ext} \, (M,L) = H^1 (C , {\mathcal O} (2A))
$$
has dimension two, when $H^0 (C,M)$ and $H^1 (C,L)$ have dimension one. Therefore there exists a nontrivial extension
$$
0 \to L \to E \to M \to 0
$$
such that $\partial$ vanishes. On the other hand, if $N \subset E$ has degree $\deg (N) > \frac{\deg (E)}{2} = 0$,
it must be contained in $M$. Since $\deg (M) = 1$ we get $N = M$ and  the extension has to be trivial.

\subsubsection{ \ } Another example, where $L = {\mathcal O}_C$ is the trivial line bundle and $M $ is the sheaf $\omega =\Omega_C^1$ of differentials on $C$, was proposed by J.~Harris (I thank B.~Mazur for explaining this to me). Choose a sextic $C' \subset {\mathbb P}^2$ with exactly two nodal singularities, and let $C$ be the normalization of $C'$. On this curve $C$ of genus 8 let $N$ be the pull-back of ${\mathcal O} (1)$ from ${\mathbb P}^2$ to $C$. One can show that there exists an extension
$$
0 \to N \to E \to \omega \, N^{-1} \to 0
$$ 
such that the boundary map
$$
H^0 (C , \omega \, N^{-1}) \to H^1 (C,N)
$$
has rank three, when $H^1 (C,N)$ has dimension four. Furthermore
 $E$ is stable and has a nowhere vanishing section. Therefore $E$ is an extension
$$
0 \to {\mathcal O}_C \to E \to \omega \to 0
$$
with associated boundary map
$$
\partial : H^0 (C,\omega) \to H^1 (C , {\mathcal O}_C)
$$
which is neither injective nor surjective.

\section{Projective subspaces in secant varieties}

\subsection{ \ }

Let $k$ be a field of characteristic zero, $C$ a smooth projective curve over $k$, and $C_{\bar k} = C \underset{k}{\otimes} \bar k$ its extension of scalars to the algebraic closure of $k$. We assume that $C_{\bar k}$ is irreducible of genus $g \geq 0$.

\medskip

Define $\Phi \subset {\mathbb Z}$ to be the image of the degree map
$$
\deg : {\rm Pic} \, (C) \to {\mathbb Z} \, .
$$
Consider a line bundle $N$ on $C$. Each cohomology class
$$
e \in H^1 (C , N^{-1}) = {\rm Ext} \, (N , {\mathcal O}_C)
$$
classifies an extension
$$
0 \to  {\mathcal O}_C \to E \to N \to 0 \, .
$$
Let $n = \deg (N)$ be the degree of $N$.

\bigskip

\noindent {\bf Proposition 2.} Assume that $n \geq 0$, that $N$ is not trivial
 and that $n + 2g-2$ lies in $2 \, \Phi$.
Let $V \subset H^1 (C,N^{-1})$ be a $k$-vector space of dimension
$$
\dim (V) \geq n-m+g \, ,
$$
where $m$ is the  integer defined by formula (3) below.

\medskip

Then there exists $e \in V$ such that the corresponding vector bundle $E$ is semi-stable (over $C_{\bar k})$.

\subsection{Proof.} Since $n+2g-2 \in 2 \, \Phi$, we can choose a line bundle $H$ on $C$ such that, if $\omega = \Omega_C^1$,
$$
\deg (N \, \omega) = 2 \deg (H) \, .
$$
If $H' = N \, \omega \, H^{-1}$, we get $\deg (H') = \deg (H)$, and
$$
N \, \omega = H \, H' \, .
$$
Since ${\rm Pic}^0 (C_{\bar k})$ is divisible, there exists a line bundle $A$ of degree zero on $C_{\bar k}$ such that
$$
H' = H \, A^2 \, .
$$
Let $M = H \, A$ (on $C_{\bar k}$) and $L = \omega \, M^{-1}$. We get
$$
N = H \, H' \, \omega^{-1} = H^2 \, A^2 \, \omega^{-1} = M (\omega \, M^{-1})^{-1} = M \, L^{-1} \, .
$$
Any class $e \in H^1 (C,N^{-1})$ defines an extension
$$
0 \to {\mathcal O}_C \to E \to N \to 0
$$
over $C$ and, by tensoring by $L$, an extension
$$
0 \to L \to E \otimes L \to M \to 0 \leqno (2)
$$
over $C_{\bar k}$. The vector bundle $E$ is semi-stable if and only if $E \otimes L$ is semi-stable.

\medskip

From now on, and till the end of \S~2, we 
assume that $k = \bar k$. Since $\deg (N) \geq 0$ we have $\deg (L) \leq \deg (M)$. Furthermore
$$
\deg (L) + \deg (M) = 2g-2 \, .
$$
Therefore, by Proposition 1, $E \otimes L$ is semi-stable if and only if the boundary map
$$
\partial_{e} : H^0 (C,M) \to H^1 (C,L)
$$
defined by (2) is an isomorphism. Note that, by Serre duality,
$$
H^1 (C,L) = H^0 (C , \omega \, L^{-1})^* = H^0 (C,M)^*
$$
has the same dimension as $H^0 (C,M)$. Let 
$$
m = \dim_k H^0 (C,M) \, . \leqno (3)
$$

\medskip

To prove Proposition 2, we now follow an argument of C.Voisin.
The map $\partial_e$ is the cup-product  by $e \in H^1 (C , L \, M^{-1})$. Therefore, by Serre duality again, the map
$$
H^1 (C , N^{-1}) = H^0 (C , N \, \omega)^* = H^0 (C , M^2)^* \to {\rm Hom} \, (H^0 (C,M) \to H^0 (C,M)^*)
$$
which maps $e$ to $\partial_e$ is dual to the cup-product
$$
H^0 (C,M)^{\otimes 2} \to H^0 (C , M^2) \, .
$$
We denote by
$$
\mu : H^0 (C,M)^{\otimes 2} \to V^* 
$$
the composite of this cup-product with the projection of $H^0 (C , M^2)$
onto the dual of $V$. Since the cup-product is commutative, any element
in $V$ defines, via $\mu$, a quadric in the projective space $\textbf{P}(H^0 (C,M))$.

Arguing by contradiction, we assume that all these quadrics are singular.
Consider the Zariski closure
$B \subset \textbf{P}(H^0 (C,M))$ of the union of the singular loci of the quadrics with  singular locus of minimal
dimension, and let $b$ be the dimension of $B$.

Let $\sigma \in H^0 (C,M)$ be a representative of a generic point $[\sigma] \in B$.
We claim that the map
$$ \mu_{\sigma} : H^0 (C,M) \to V^* $$
mapping $\tau \in H^0 (C,M)$ to 
$$ \mu_{\sigma}(\tau) =  \mu (\sigma \otimes \tau) $$
has rank at most $b$. Indeed, it follows from the definitions that a
quadric $q \in V$ is singular at $\tau \in H^0 (C,M)$ if and only if
it lies in the  subspace $Q_{\tau} \subset V$ orthogonal to the image of
$\mu_{\tau}$. Generically, the singular locus of $q$ is minimal. 
Therefore the union all the vector spaces $Q_{\tau}$, $[\tau] \in B$,
is an open subset of $V$. Since $[\sigma]$ is generic in $B$,
the dimension of $Q_{\sigma}$ is at least ${\rm {dim}}(V) - b$,
and the rank of $\mu_{\sigma}$ is at most $b$ as claimed.

This implies  that the kernel $H_{\sigma}\subset H^0 (C,M)$ of $\mu_{\sigma}$
has dimension $c \geq m - b$ (note that
this dimension $c$ has a fixed value when $[\sigma]$ is generic
in $B$).
Let $K \subset H^0 (C , M^2)$ be the  subspace orthogonal to $V$. 
By definition, the vector space
$$K_{\sigma} = \sigma \cup  H_{\sigma}$$
is contained in $K$. Its dimension is $c$.

\medskip

On the other hand, we can choose points $x_0, ..., x_b$ on $C$ and vectors
 $\sigma_0, ... , \sigma_b \in H^0 (C,M)$ such that $[\sigma_i]$ lies in
 $B$ and 
 $$\sigma_i(x_j) = \delta_{ij}$$ 
 for all $i$ and $j$. By moving $x_i$ without moving the other points,
 we can also assume that, for every $i$, at least one section in
 $H_{\sigma_i}$ does not vanish at $x_i$. As a consequence,
 $K_{\sigma_i}$ is not contain in the sum of the $K_{\sigma_j}$'s,
 $j \neq i$, and the dimension of the sum of the subspaces $K_{\sigma_i}$,
 $i = 0, ... ,b$, is at least 
 $$ b + c \geq m \, . $$
 Therefore $K$ has dimension at least $m$ and, since $H^1 (C , N^{-1})$
 has dimension $n + g - 1$, the dimension of $V$ is at most $n + g - m - 1$.
 This contradicts our hypotheses.

\subsection{Remarks.}

\subsubsection{ \ } In the proof of Proposition 2, $N = M^2 \, \omega^{-1}$, therefore
$$
\deg (M) = \frac{n}{2} + g - 1 \geq g-1 \, .
$$
By the Riemann-Roch theorem:
$$
\chi (C,M) = \deg (M) - g + 1 = \frac{n}{2} \, .
$$
By Clifford's theorem
$$
\dim_k H^1 (C,M) \leq {\rm Sup} (g-1 , 0) \, ,
$$
and $\dim_k H^1 (C,M) = 0$ whenever $\deg (M) > 2g-2$. Therefore
$$
\frac{n}{2} \leq m \leq \frac{n}{2} + {\rm Sup} (g-1,0) \, ,
$$
hence
$$
\frac{n}{2} + {\rm Inf} (g,1) \leq n-m+g \leq \frac{n}{2} + g
$$
and $n-m+g = \frac{n}{2} + g$ as soon as $n > 2g-2$.

\subsubsection{ \ } If $C$ has a $k$-rational point $\Phi = {\mathbb Z}$. In that case, the conditions on $n$ in Proposition 2 mean that $n$ is a nonnegative even integer.

\subsection{Secant varieties}

\subsubsection{ \ } The Proposition 2 can be rephrased in terms of secant varieties. Let $d \geq 1$ be an integer and
$$
n = 2 \, d + 2 \, .
$$
Let $k$ be an algebraically closed field of characteristic zero and $C$ a smooth connected projective curve over $k$,
of genus $g$, say $C \subset {\mathbb P}$. Let $N^{-1} = \omega \, {\mathcal O} (-1)$ be the Serre dual of the canonical sheaf on $C$, and assume that $\deg (N) = n$.

\medskip

Consider the secant variety
$$
\Sigma_d = \bigcup_{Z \in X^{(d)}} \langle Z \rangle \, ,
$$
swept out by the linear spans of $d$-uples of points on $C$. Define $m$ as in (3).

\bigskip

\noindent {\bf Theorem 1.} The secant variety $\Sigma_d$ does not contain any projective space ${\mathbb P}^{\delta}$ of dimension
$$
\delta \geq n-m+g-1 \, .
$$

\subsubsection{Proof.} Let $e \in H^1 (C,N^{-1})$, $e \ne 0$, and
$$
0 \to {\mathcal O}_C \to E \to N \to 0
$$
the corresponding extension. the semi-stability of  $E$  means that $e$ does not lie in the image of the boundary map
$$
\partial_D : H^0 (D , N^{-1} (D)) \to H^1 (C,N^{-1})
$$
coming from
$$
0 \to N^{-1} \to N^{-1} (D) \to N^{-1} (D) / D \to 0 \, ,
$$
for any effective divisor $D$ of degree less than $\frac{n}{2}$, {\it i.e.} $\deg (D) \leq d$. This condition happens to be equivalent to the fact that the point in ${\mathbb P} = {\mathbb P} (H^1 (C , N^{-1}))$ defined by $e$ does not belong to $\Sigma_d$. For more details see \cite{B}, p.~451, or \cite{S2}, \S~1.6. Therefore Theorem 1 follows from Proposition 2.

\subsubsection{ \ } Using 2.4.1 we see that the lower bound
$$
\delta_0 = n-m+g-1
$$
in Theorem 1 is such that
$$
\delta_0 \geq \frac{n}{2} + {\rm Inf} \, (g-1,0) = d + {\rm Inf} \, (g,1) \, ,
$$
and $\delta_0 = d+g$ when $n > 2g-2$. The remark 1.3 above suggests that this bound is not optimal. According to C.~Voisin, when $g > 0$, Theorem 1 should remain true with $\delta \geq d$ (\cite{S2} , \S~1.3).

\section{Semi-stable extensions on arithmetic surfaces}

\subsection{ \ } Let $F$ be a number field, ${\mathcal O}_F$ its ring of integers and $S = {\rm Spec} ({\mathcal O}_F)$. Consider a semi-stable curve $X$ over $S$ such that $X$ is regular and its generic fiber $X_F$ is geometrically irreducible of genus $g \geq 0$. Let
$$
\deg : {\rm Pic} \, (X) \to {\mathbb Z}
$$
be the morphism sending the class of a line bundle on $X$ to the degree of its restriction to $X_F$. Call $\Phi \subset {\mathbb Z}$ the image of that degree map.

\medskip

Let $\bar N = (N,h)$ be an hermitian line bundle over $X$, {\it i.e.} a line bundle $N$ on $X$ together with an hermitian metric $h$ on the restriction $N_{\mathbb C}$ of $N$ to $X({\mathbb C})$ which is invariant under complex conjugation. The cohomology group
$$
\Lambda = H^1 (X , N^{-1})
$$
is a finitely generated module over ${\mathcal O}_F$. For every complex embedding $\sigma : F \to {\mathbb C}$, let $X_{\sigma} = X \otimes {\mathbb C}$ be the corresponding surface and $\Lambda_{\sigma} = \Lambda \otimes {\mathbb C}$. This cohomology group
$$
\Lambda_{\sigma} = H^1 (X_{\sigma} , N_{\mathbb C}^{-1}) 
$$
is canonically isomorphic to the complex vector space $\Omega^1 (X_{\sigma} , N_{\mathbb C}^{-1})$ of holomorphic differential froms
 with coefficients in the restriction $N_{\mathbb C}^{-1}$ of the line bundle $N^{-1}$ to $X ({\mathbb C}) = \underset{\sigma}{\coprod} \, X_{\sigma}$. Given $\alpha \in \Omega^1 (X_{\sigma} , N_{\mathbb C}^{-1})$, we let $\alpha^*$ be its transposed conjugate (the definition of which uses the metric $h$), and we define
$$
\Vert \alpha \Vert_{L^2}^2 = \frac{i}{2\pi} \int_{X_{\sigma}} \alpha^* \, \alpha \, .
$$
Given $e \in \Lambda$, we let
$$
\Vert e \Vert = \underset{\sigma}{\rm Sup} \, \Vert \sigma (e) \Vert_{L^2} \, ,
$$
where $\sigma$ runs over all complex embeddings of $F$.

\medskip

We are interested in (the logarithm of) the successive minima of $\Lambda$. Namely, for any positive integer $k \leq rk(\Lambda)$, we let $\mu_k$ be the infimum of all real numbers $\mu$ such that there exist $k$ elements $e_1 , \ldots , e_k$ in $\Lambda$ which are linearly independent in
$$
\Lambda \otimes F = H^1 (X_F , N^{-1})
$$ 
and such that
$$
\Vert e_i \Vert \leq \exp (\mu) \qquad \hbox{for all} \ i = 1,\ldots , k \, . \leqno (4)
$$

Let $n = \deg (N)$. We assume that $n > 0$ and that $n+2g-2 \in 2 \, \Phi$. We define $m$ by the formula (3) above (with ground field $F$ instead of $k$). Finally, let
$$
d = [F:{\mathbb Q}]
$$
be the degree of $F$ over ${\mathbb Q}$.

\bigskip

\noindent {\bf Theorem 2.} Assume that $k \geq n-m+g$. Then 
$$
\mu_k \geq \frac{\hat c_1 (\bar N)^2}{2 \, n \, d} - A \, ,
$$
where
$$
A = \frac{1}{n \, d} + \log (m (n + g - 1)) \, ,
$$
and $\hat c_1 (\bar N)^2 \in {\mathbb R}$ denotes the self-intersection of the arithmetic Chern class $\hat c_1 (\bar N) \in \widehat{\rm CH}^1 (X)$.

\subsection{Proof.} Let $e_1 , \ldots , e_k$ be elements of $\Lambda$ which are $F$-linearly independent and such that (4) holds. Call $V \subset H^1 (X_F , N^{-1})$ the $F$-vector space spanned by $e_1 , \ldots , e_k$. According to Proposition 2 there exists $e \in V$ such that the corresponding extension $E$ of $N$ by the trivial line bundle on $X_F$ is semi-stable on $X_{\bar F}$. Furthermore, using the notation of the proof of Proposition 2, $E$ is semi-stable as soon as
$$
\partial_e : H^0 (X_{\bar F} , M) \to H^1 (X_{\bar F} , L)
$$
is an isomorphism. Choosing a basis of these two vector spaces, we get a polynomial $P$ of degree $m$ on $V \otimes \bar F$ such that
$$
P(e) = \det (\partial_e) \, ,
$$
so that $E$ is semi-stable as soon as $P(e) \ne 0$. Therefore, by a standard argument (see \cite{S2}, proof of Proposition 5), there exists $k$ integers $n_1 , \ldots , n_k$, with $\vert n_i \vert \leq m$ for all $i$, such that
$$
e = n_1 \, e_1 + \ldots + n_k \, e_k \leqno (5)
$$
satisfies $P(e) \ne 0$, hence $E$ is semi-stable.

\medskip

From the definition of $\mu_k$ and (5) we get
$$
\Vert e \Vert \leq m \, k \exp (\mu_k) \leq m (n+g-1) \exp (\mu_k) \leqno (6)
$$
(since $rk (\Lambda) = n+g-1$). According to a result of Miyaoka (\cite{M2}, \cite{S1} Theorem 1) and Moriwaki \cite{Mo}, this implies that, for any choice of a metric on $E$ (invariant under complex conjugation), the inequality ``\`a la Bogomolov''
$$
\hat c_1 (\bar E)^2 \leq 4 \, \hat c_2 (\bar E) \leqno (7)
$$
is satisfied in ${\mathbb R}$. Here, as in \cite{S1} \S~2.1, given $x \in \widehat{\rm CH}^2 (X)$ we also denote by $x \in {\mathbb R}$ its arithmetic degree $\widehat\deg \, (x)$.

\medskip

We now proceed in a way similar to \cite{S1}, Proposition 1 and Corollary (where more details can be found). Recall that $E$ is an extension
$$
0 \to {\mathcal O}_X \to E \to N \to 0 \, . \leqno (8)
$$
We endow ${\mathcal O}_X$ with the trivial metric and $N$ with a metric $h'$ to be specified below.
 For any choice of a smooth splitting of (8) over $X({\mathbb C})$, we get a metric on $E$, namely the orthogonal direct sum of the chosen metrics on ${\mathcal O}_X$ and $N$. The Cauchy-Riemann operator on $E_{\mathbb C}$ can be written in matrix form according to that splitting:
$$
\bar\partial_E = \begin{pmatrix} \bar\partial &\alpha \\ 0 &\bar\partial_N \end{pmatrix} \, ,
$$
where $\alpha$ is a smooth form of type $(0,1)$ over $X({\mathbb C})$ with coefficients in $N_{\mathbb C}^{-1}$. One can choose the smooth splitting of (8) over $X({\mathbb C})$ in such a way that $\alpha$ is the harmonic representative of the restriction of $e$ to $X({\mathbb C})$. With this choice we get
$$
\hat c_1 (\bar E) = \hat c_1 (N,h')
$$
and
$$
2 \, \hat c_2 (\bar E) = \sum_{\sigma : F \to {\mathbb C}} \Vert \sigma (e) \Vert'^2_{L^2} \, ,
$$
where, for every complex embedding $\sigma$, $\Vert \cdot \Vert'_{L^2}$ is the $L^2$-norm on $H^1 (X_{\sigma} , N_{\mathbb C}^{-1})$ defined by $h'$.

\medskip

Now let $t = \Vert e \Vert^2$ and let us choose $h' = th$. We get
$$
\hat c_1 (N,h')^2 = \hat c_1 (\bar N)^2 - n \, d \log (t)
$$
and
$$
\sum_{\sigma} \Vert \sigma (e) \Vert'^2_{L^2} = t^{-1} \sum_{\sigma} \Vert \sigma (e) \Vert_{L^2}^2 \leq 1 \, .
$$
Therefore the inequality (7) reads
$$
\hat c_1 (\bar N)^2 \leq 2 \, n \, d \log \Vert e \Vert + 2 \, .
$$
Since, by (6),
$$
\log \Vert e \Vert \leq \mu_k + \log (m (n+g-1)) \, ,
$$
Theorem 2 follows.

 {\it CNRS and IH\'ES}

 \bigskip

 35 Route de Chartres

91440  Bures sur Yvette, France

\bigskip

soule@ihes.fr

\bigskip

\newpage

\begin{thebibliography}{999}
\bibitem{ACGH} E.Arbarello, M.Cornalba, P.A.Griffiths, J.Harris:
Geometry of Algebraic Curves, Vol.~I (1985), Springer-Verlag.

\bibitem{B} A.Bertram: Moduli of rank $2$ vector bundles, theta
divisors, and the geometry of curves in projective space, {\it J.
Diff. Geom.} {\bf 35} (1992), 429-469.


\bibitem{M1} Y.Miyaoka: Talk in MPI Bonn, March 1988.

\bibitem{M2} Y.Miyaoka: Bogomolov inequality on arithmetic surfaces,
talk at the Oberwolfach conference on ``Arithmetical Algebraic
Geometry'', G. Harder and N. Katz org. (1988).


\bibitem{MB} L.Moret-Bailly:
Hauteurs et classes de Chern sur les surfaces arithm\'etiques.
Les pinceaux de courbes elliptiques, S\'emin., Paris 1988, L.Szpiro org.,
 Ast\'erisque \textbf{183}  (1990), 37-58.

\bibitem{Mo} A.Moriwaki:
Inequality of Bogomolov-Gieseker type on arithmetic surfaces.
{\it Duke Math. J.} \textbf{74}, No.3 (1994), 713-761.

\bibitem{P} A.N.Parshin:
The Bogomolov-Miyaoka-Yau inequality for the arithmetical surfaces and its applications.
S\'emin. Th\'eor. Nombres, Paris 1986-87, Prog. Math. 75 (1988), 299-312.

\bibitem{S1} C.Soul\'e: A vanishing theorem on arithmetic surfaces,
{\it Invent. Math.} {\bf 116}  (1994), 577-599.

\bibitem{S2} C.Soul\'e: Secant varieties and successive minima,
 {\it J. Algebraic Geom.} \textbf{13},  no. 2 (2004), 323--341.

\bibitem{V} C.Voisin:
On linear subspaces contained in the secant varieties of a projective curve,
{\it J. Algebraic Geom.}  \textbf{13 },  no. 2 (2004), 343--347.

\end{thebibliography}
\end{document}